%
% file:         tams.tex (former tran561.tex)
% created:      received from AMS dec 12 2008
% modified:     pasha jan 22 2014
% modification: fixed a typo; removed address, put current email
%

\input amstex
\documentstyle{tran}
\NoBlackBoxes

\pageheight{24.1cm}

\loadeusm
\def\scr#1{{\fam\eusmfam\relax#1}}
\def\germ{\frak}

\issueinfo{334}
  {1}
  {November}
  {1992}

\pageno=143

\topmatter

\title Central extensions of current algebras\endtitle
\author Paul Zusmanovich\endauthor
\address  \endaddress
%\curraddr Department of Mathematics, Bar-Ilan University,
%Ramat-Gan 52900, Israel \endcurraddr
\email justpasha\@gmail.com \endemail
\date April 18, 1990 and, in revised form, September 4, 1990. Last minor \
revision for arXiv January 22, 2014 \enddate
\abstract The second cohomology group of Lie algebras of kind $L\otimes U$ with
trivial coefficients is investigated, where $L$ admits a decomposition with
one-dimensional root spaces and $U$ is an arbitrary associative commutative
algebra with unit. This paper gives a unification of some recent results of C.
Kassel and A. Haddi and provides a determination of central extensions of
certain modular semisimple Lie algebras.
\endabstract
\endtopmatter

\document
\heading Introduction\endheading
\par The main result of this paper (Theorem 1.3) gives the expression of the
second cohomology group of a Lie algebra $L\otimes U$ with trivial
coefficients (denoted by $H^2(L\otimes U))$, where $L$ belongs to the class of
Lie algebras referred to as Block algebras (for the definition, see below), and
$U$ is a commutative associative algebra with unit, through the first-order
cyclic cohomology $HC^1(U)$ and the dual $U^*$. Since the classical Lie
algebras, Zassenhaus algebras $W_1(n)$, and infinite-dimensional Witt algebras
$W$ are Block algebras, we obtain as corollaries the values of $H^2(L\otimes U)$
in all these cases. In the second section of this paper the cohomology groups
$H^2(W_1(n)\otimes U)$ are computed by a different method, based on the
presentation of $W_1(n)$ as a deformation of tensor product of $W_1(1)$ and
divided power algebra. Notice that $W_1(\infty)$ is not a Block algebra and
hence this case is not covered by the previous theorem. The third section
contains an application of the preceding results to the computation of central
extensions of certain modular semisimple Lie algebras.
\par Note that earlier C. Kassel \cite9 and A. Haddi \cite8 computed the second
homology group $H_2(L\otimes U)$, where $L$ is classical or $W$ and its
subalgebra $W_1$ respectively. Their method is quite different from ours and is
based on the idea of universal central extensions. Our main theorem may be
considered as a generalization and unification of Kassel's and Haddi's results.
\par All algebras are defined over a field $\Bbb F$ of characteristic
$p\not=2,3$, if the otherwise condition is not stated. Remember that $U$
denotes an arbitrary associative commutative algebra with a unit 1.
\heading1. The main theorem and corollaries\endheading
\par Let $T$ be a torus in a Lie algebra $L$. Assume that there exists a root
space decomposition of $L$ with respect to the $T$-action: %\longpage{13pc}
\pagebreak
$$L=T\oplus\bigoplus_{\alpha \in R}L_\alpha .\tag1$$
Assume, moreover, that all root spaces are one-dimensional: $L_\alpha
=\Bbb F e_\alpha$ and \linebreak $\alpha ([L_\alpha ,L_{-\alpha }])\not=0$ for all $\alpha
\in R$. Set $h_\alpha =[e_{-\alpha },e_\alpha ]$, $[e_\alpha ,e_\beta
]=N_{\alpha ,\beta }e_{\alpha +\beta }$, $N_{\alpha ,\beta }\in\Bbb F$. The
next assumption concerning the structure of $L$ is $[L,L]=L$. This is
equivalent to $T=\sum\Bbb F h_\alpha $. Let $\{h_\alpha |\alpha \in B\}$,
$B\subseteq R$, be a basis in $T$.
\par Lie algebras $L$ satisfying all these properties will be called {\it Block
algebras\/}. In \cite1 R. E. Block showed that all finite-dimensional Block
algebras over a perfect field of characteristic $p>5$ with $Z(L)=0$ are direct
sums of simple classical Lie algebras and Albert-Zassenhaus algebras.
\par In the sequel, $L$ will denote a Block algebra except in Lemma 1.1 and
Lemma 3.1, where $L$ is arbitrary.
\par Our cohomological notions and notations are standard (cf.\ \cite7 for Lie
algebra cohomology and \cite5 for cyclic cohomology). When we consider the Lie
algebra cohomology, we use the standard cochain complex. Notice that since $U$
is commutative, the first-order cyclic cohomology $HC^1(U)$ is exactly the
module of all skew-symmetric functions $F\in(U\otimes U)^*$ such that
$$F(uv,w)+F(wu,v)+F(vw,u)=0.$$
If $\scr D$ is a subalgebra in $\operatorname{Der}(U)$, the module of $\scr
D$-invariant cyclic cohomology $HC^1(U)^{\scr D}$ consists of cocycles
satisfying
$$F(D(u),v)+F(u,D(v))=0$$
for all $D\in \scr D$.
\par Recall that the Lie algebra structure on $L\otimes U$ is defined by
$[x\otimes u,y\otimes v]=[x,y]\otimes uv$.
\par Now compute $H^2(L\otimes U)$.
\par The following lemma is probably one of the most useful results
for computation of Lie algebra cohomology groups.
\proclaim{Lemma 1.1. \rm(Cf\. \cite{7, Theorem 1.5.2})}
 Let $L$ be a Lie algebra,
and let $T$ be a torus in $L$ such that there exists a root space decomposition
$(1)$. Then for each cohomology class $[\Phi]\in H^k(L)$ there exists 
a representative cocycle $\Phi$ such that $\Phi(L_{\alpha _1},\dotsc,
L_{\alpha _k})=0$ if $\alpha _1+\dotsb+\alpha _k\not=0$.\qed
\endproclaim
\proclaim{Lemma 1.2 \rm$H^2(\operatorname{sl}(2)\otimes U)\cong HC^1(U)$}
\par If we choose a basis $\{e,f,h\}$ in $\operatorname{sl}(2)$ with
multiplication table 
$$[h,e]=-\alpha e,\qquad [h,f]=\alpha f,\qquad [e,f]=h$$
for certain $\alpha \in\Bbb F\backslash\{0\}$, then the basic cocycles may be
chosen as follows{\rm:}
$$\gathered
{\aligned
h\otimes u\wedge h\otimes v\mapsto&-\alpha F(u,v),\\
e\otimes u\wedge f\otimes v\mapsto& F(u,v),\endaligned}\\
\text{other combinations }\mapsto 0\endgathered$$
where $F\in HC^1(U)$.
\endproclaim
\demo{Proof} Straightforward computations or reference to Kassel \cite9.
Perhaps it should only be remarked that the coboundary cyclic $(=$ Hochschild)
cohomology condition arises from the Jacobi identity.\qed
\enddemo
\pagebreak
\proclaim{Theorem 1.3} $H^2(L\otimes U)\cong (\bigoplus_{\alpha \in
R}HC^1(U)_\alpha )/I\oplus(\bigoplus_{\alpha \in R}U^*_\alpha )/J$ where 
$HC^1(U)_\alpha \mathbreak
\cong HC^1(U)$, $U^*_\alpha \cong U^*$ for any $\alpha \in R$,
$I$ is the linear span of the following elements{\rm:}
$$\alpha (h)F_\alpha -\sum_{\beta \in B}\lambda _{\alpha \beta }\beta
(h)F_\beta ,\quad
\alpha (h_\beta )F_\alpha -\beta (h_\alpha )F_\beta ,\quad N_{\beta ,-\alpha
-\beta }F_\alpha +N_{\alpha ,-\alpha -\beta }F_\beta $$
for all $h_\alpha, h\in T$, where $h_\alpha =\sum_{\beta \in B}\lambda
_{\alpha \beta }h_\beta $ is the linear expression of $h_\alpha $ in terms of
the basis $\{h_\beta |\beta \in B\}$, and for all $\alpha $, $\beta \in R$,
where $F_\alpha $ is an arbitrary element in $HC^1(U)_\alpha $. $J$ is the
linear span of the elements
$$\multline
N_{\beta ,\gamma}G_\alpha +N_{\gamma,\alpha }G_\beta +N_{\alpha ,\beta }G_
\gamma,\quad \alpha ,\beta ,\gamma\in R,\ \alpha +\beta +\gamma=0;\\
G_{-\alpha }+G_\alpha , \alpha \in R; G_\alpha , \alpha \in
B,\endmultline
$$
for all $G_\alpha \in U^*_\alpha $.
\par Basic cocycles may be chosen as follows{\rm:}
$$\gather
h_\alpha \otimes u\wedge h_\beta \otimes v\mapsto -\alpha (h_\beta )F_\alpha
(u,v),\qquad \alpha ,\beta \in B,\tag2\\
e_{-\alpha }\otimes u\wedge e_\alpha \otimes v\mapsto F_\alpha (u,v)+G_\alpha
(uv),\qquad \alpha \in R,\tag3\\
\text{other combinations }\mapsto 0.\tag4
\endgather
$$
\endproclaim
\demo{Proof} Obviously for any root $\alpha $ it is possible to find a
complement $Z_\alpha $ of $\Bbb F h_\alpha $ in $T$ such that $[Z_\alpha
,e_\alpha ]=[Z_\alpha ,e_{-\alpha }]=0$. $T\oplus\Bbb F e_{-\alpha }\oplus 
\Bbb F e_\alpha $ is a subalgebra in $L$ which is isomorphic to an
(unessential) central extension of $\operatorname{sl}(2)$. Since
$$H^2((T\oplus \Bbb F e_{-\alpha }\oplus\Bbb F e_\alpha )\otimes U)\cong 
H^2(Z_\alpha \otimes U)+H^2((\Bbb F h_\alpha \oplus \Bbb F e_{-\alpha }\oplus 
\Bbb F e_\alpha )\otimes U),$$
by Lemma 1.2 we get
$$H^2((T\oplus \Bbb F e_{-\alpha }\oplus \Bbb Fe_\alpha )\otimes U)\cong
C^2(Z_\alpha \otimes U)+HC^1(U).$$
Hence for arbitrary $\Phi\in Z^2(L\otimes U)$ and $\alpha \in R$ we have 
$$\gather
\Phi(h_\alpha \otimes u,h_\alpha \otimes v)=-\alpha (h_\alpha )F_\alpha
(u,v),\\
\Phi(e_{-\alpha }\otimes u,e_\alpha \otimes v)=F_\alpha (u,v)+G_\alpha (uv),\\
\Phi(h_\alpha \otimes U,Z_\alpha \otimes U)=0,\endgather
$$
for certain $F_\alpha \in HC^1(U)$ and $G_\alpha \in U^*$.
\par For any $h\in T$, $h-(\alpha (h)/\alpha (h_\alpha ))h_\alpha \in Z_\alpha
$. Hence
$$\Phi(h_\alpha \otimes u,(h-(\alpha (h)/\alpha (h_\alpha ))h_\alpha )\otimes
v)=0$$
and
$$\Phi(h_\alpha \otimes u,h\otimes v)=(\alpha (h)/\alpha (h_\alpha
))\Phi(h_\alpha \otimes u,h_\alpha \otimes v)=-\alpha (h)F_\alpha (u,v).\tag5$$
Now suppose that in the last equality $h_\alpha =\sum_{\beta \in B}\lambda
_{\alpha \beta }h_\beta $. Then 
$$\align
-\alpha (h)F_\alpha (u,v)=&\sum_{\beta \in B}\lambda _{\alpha \beta
}\Phi(h_\beta \otimes u,h\otimes v)\\
=&-\sum_{\beta \in B}\lambda _{\alpha \beta }\beta (h)F_\beta (u,v).\endalign
$$
This implies the defining relations for $I$ of the first kind. Further,
substituting in (5) $h=h_\beta $, we obtain (2) and (because of the
skew-symmetry of $\Phi)$ the defining relations for $I$ of the second kind:
$$\alpha (h_\beta )F_\alpha (u,v)=\beta (h_\alpha )F_\beta (u,v).$$
The cocycle equation $d\Phi=0$ for the triple $e_\alpha \otimes u$, $e_\beta
\otimes v$, $e_\gamma\otimes w$, $\alpha +\beta +\gamma=0$, is equivalent to
the relation $N_{\alpha ,\beta }F_\gamma(uv,w)+N_{\alpha ,\beta
}G_\gamma(uvw)+$(simultaneous cyclic permutations of $(\alpha ,\beta ,\gamma)$
and $(u,v,w))=0$. Substituting here $v=w=1$, we obtain
$N_{\beta ,\gamma}G_\alpha (u)+N_{\gamma,\alpha }G_\beta (u)+N_{\alpha ,\beta
}G_\gamma(u)=0$ (which is equivalent to the defining relations for $J$ of the
first kind) and hence
$$N_{\beta ,\gamma}F_\alpha (vw,u)+N_{\gamma,\alpha }F_\beta (wu,v)+N_{\alpha
,\beta }F_\gamma (uv,w)=0.$$
Substituting in the last equality $w=1$ we get the defining relations for $I$
of the third kind:
$$N_{\beta ,-\alpha -\beta }F_\alpha (v,u)+N_{-\alpha -\beta ,\alpha }F_\beta
(u,v)=0.$$
Since $T\otimes 1$ is a torus in $L\otimes U$, by Lemma 1.1 we may assume
$$\gather
\Phi(e_\alpha \otimes U,e_\beta \otimes U)=0\quad\text{if }\alpha +\beta
\not=0,\\
\Phi(T\otimes U,e_\alpha \otimes U)=0,\endgather
$$ 
which implies (4).
\par Now put $\omega \in C^1(L\otimes U)$ as follows:
$$\omega (h_\alpha \otimes u)=G_\alpha (u),\quad\alpha \in B;\qquad \omega
(e_\alpha \otimes U)=0,\quad \alpha \in R.$$
Replacing $\Phi$ by $\Phi-d\omega $, we obtain a cocycle which satisfies all
previous equalities and, moreover, $G_\alpha =0$, $\alpha \in B$. Finally, the
relation $G_\alpha +G_{-\alpha }=0$ obviously follows from the skew-symmetry of
$\Phi$.
\par It remains to prove only that $\Phi$ described in the theorem is really a
cocycle and that all such cocycles are independent. This is achieved by an easy
argument which mainly repeats the previous one and is left to the reader.
\qed
\enddemo
\par We keep the symbol $\germ g$ for classical Lie algebras, that is, for Lie
algebras obtained from integral forms of complex classical simple Lie algebras
by tensoring on $\Bbb F$. Let $R$ be a root system of $\germ g$, $B$ be a
basis, $\{h_\alpha ,\ \alpha \in B;\ e_\alpha ,\ \alpha \in R\}$ be a Chevalley
basis of $\germ g$, and let $\langle \cdot,\cdot\rangle $ denotes a
scalar product in $R$.
\proclaim{Corollary 1.4 \cite9} $H^2(\germ g\otimes U)\cong HC^1(U)$.
\par Basic cocycles may be chosen as follows{\rm:}
$$\gathered
h_\alpha \otimes u\wedge h_\beta \otimes v\mapsto-\frac{2\langle \alpha ,\beta
\rangle }{\langle \alpha ,\alpha \rangle \langle \beta ,\beta \rangle
}F(u,v),\qquad \alpha ,\beta \in B,\\
e_{-\alpha }\otimes u\wedge e_\alpha \otimes v\mapsto \frac1{\langle \alpha
,\alpha \rangle }F(u,v),\qquad \alpha \in R,\\
\text{other combinations}\mapsto0,\endgathered\tag6
$$
where $F\in HC^1(U)$.
\endproclaim
\demo{Proof} There is the following property of the structure constants:
$$\frac{N_{\alpha ,\beta }}{\langle \gamma,\gamma\rangle }=\frac{N_{\beta
,\gamma}}{\langle \alpha ,\alpha \rangle }=\frac{N_{\gamma,\alpha }}{\langle
\beta ,\beta \rangle }\quad\text{if }\alpha +\beta +\gamma=0\tag7$$
(cf.\ \cite{4, Chapter VIII, \S2, Exercise 4}). This property and the defining
relations for $J$ of the first kind implies 
$$G_{-\alpha -\beta }=-\frac{\langle \alpha ,\alpha \rangle }{\langle \alpha
+\beta ,\alpha +\beta \rangle }G_\alpha -\frac{\langle \beta ,\beta \rangle }{
\langle \alpha +\beta ,\alpha +\beta \rangle }G_\beta ,\qquad \alpha ,\beta
,\alpha +\beta \in R. \pagebreak$$
But $G_\alpha =0$ for all $\alpha \in B$ and, consequently, $G_\alpha =0$ for
all negative roots $\alpha $. The defining relations for $J$ of the second kind
yield $G_{-\alpha }=-G_\alpha $ and hence $G_\alpha =0$ for all positive roots
$\alpha $.
\par Further, (7) together with the defining relations for $I$ of the third
kind gives $\langle \alpha ,\alpha \rangle F_\alpha =\langle \beta ,\beta
\rangle F_\beta $, $\alpha ,\beta \in R$. Hence
$$F_\alpha \frac1{\langle \alpha ,\alpha \rangle }F,\qquad \alpha \in R,\tag8$$
for certain $F\in HC^1(U)$. Since for $\alpha $, $\beta \in R$, $\alpha
(h_\beta )=2\langle \alpha ,\beta \rangle /\langle \beta ,\beta \rangle $, the
right-hand side of formula (2) takes the form
$$-\frac{2\langle \alpha ,\beta \rangle }{\langle \alpha ,\alpha \rangle
\langle \beta ,\beta \rangle }F(u,v)$$
which agrees with (6). Finally, it is easy to check that under condition (8)
the defining relations for $I$ of the first and second kinds always holds.
\qed
\enddemo
\par Now recall the definition of the modular Zassenhaus algebra $W_1(n)$. It
may be defined in two ways. Firstly as an algebra with basis $\{f_\alpha
|\alpha \in\Bbb G \Bbb F (p^n)\}$ and Lie bracket given by
$$[f_\alpha ,f_\beta ]=(\beta -\alpha )f_{\alpha +\beta }.\tag9$$
Secondly, as an algebra with basis $\{e_i|-1\le i\le p^n-2\}$ and bracket 
$$[e_i,e_j]=\left(\left(\matrix i+j+1\\j\endmatrix\right)-\left(\matrix
i+j+1\\ i\endmatrix \right)\right)e_{i+j}.\tag10$$
There are some related infinite-dimensional algebras. If we assume that the
indices in formula (9) run over all integers, then we obtain the
infinite-dimensional Witt algebra $W$ (considered over a field of
characteristic zero). If we assume that the indices in formula (10) run from
$-1$ to infinity, we obtain the infinite-dimensional algebra $W_1(\infty)$
(considered over a field of positive characteristic). Obviously $W_1(n)$,
$n<\infty$, and $W$ are Block algebras.
\proclaim{Corollary 1.5} $(\roman i)$ $(p>3,\ n<\infty)$ $H^2(W_1(n)\otimes
U)\cong H^2(W_1(n))\otimes U^*$.
\par$(\roman{ii})$ \cite8 $(p=0)$ $H^2(W\otimes U)\cong H^2(W)\otimes U^*$.
\endproclaim
\demo{Proof} (i). Choose a basis $\{f_\alpha |\alpha \in\Bbb G \Bbb F (p^n)\}$.
We have (in terms of Block algebras): $h_\alpha =2f_0$ for all $\alpha \in R=
\Bbb G \Bbb F (p^n)\backslash \{0\}$, $\alpha (f_0)=\alpha $, $N_{\alpha ,\beta
}=\beta -\alpha $. Apply Theorem 1.3. The defining relations for $I$ of the
second and third kinds give $\alpha F_\alpha =\beta F_\beta $ and $(\alpha
+2\beta )F_\alpha +(2\alpha +\beta )F_\beta =0$. These equalities evidently
imply $F_\alpha =0$ for all $\alpha $. The defining relations for $J$ give
$$(\beta -\alpha )G_{-\alpha -\beta }-(\alpha +2\beta )G_\alpha +(2\alpha
+\beta )G_\beta =0,\quad G_{-\alpha }=-G_{\alpha},\quad G_1=0.$$
Then we proceed, as in Block's paper \cite{2, p. 1449}, to obtain
$G_\alpha =(\alpha -\alpha ^3)G$ for certain $G\in U^*$. Together with a result
of Block \cite{2, Theorem 5.1} this completes the proof.
\par(ii) Quite analogously.\qed
\enddemo
\par Corollaries 1.4 and 1.5 show that the ``cyclic cohomology'' summand in our
main formula (Theorem 1.3) arises from the ``classical'' properties of Block
algebras, while the ``dual'' summand arises from ``Zassenhaus''
 $=($``nonclassical'') ones. Of course, Theorem 1.1 may be applied to compute
$H^2(L\otimes U)$ for other classes of Block algebras, e.g. Albert-Zassenhaus
or Virasoro (both modular and nonmodular).
\heading2. Central extensions of $W_1(n)\otimes U$\endheading
\par Define a Lie algebra $\germ L(U,D)$, where $D\in\operatorname{Der}(U)$, as
the vector space $W_1(1)\otimes U$ with Lie bracket given by $\{x,y\}=[x,y]
+\Psi(x,y)$ where $[\cdot,\cdot ]$ is the ordinary Lie bracket in
$W_1(1)\otimes U$ and 
$$\Psi(e_i\otimes u,e_j\otimes v)=\cases
e_{p-2}\otimes vD(u)-uD(v),\quad i=j=-1,\\
0,\quad\text{otherwise}.\endcases
$$
\par Now the notion of divided power algebra $\scr O_1(n)$ is necessary. It is
defined as an associative commutative algebra over a field of characteristic
$p>0$ with basis $\{x^i|\ 0\le i<p^n\}$ and multiplication $x^ix^j=\binom{
i+j}{j}x^{i+j}$ (the case $n=\infty$ is also included). It is known that $
\scr O_1(n)$ is isomorphic to the truncated polynomial algebra $\scr O_n=\Bbb
F[x_1,\dotsc, x_n]/(x^p_1,\dotsc, x^p_n)$. We reserve the special notation
$\partial $ for the derivation of $\scr O_1(n)$ defined by $\partial
(x^i)=x^{i-1}$.
\par The following remarkable observation is due to M. I. Kuznetsov \cite{10}:
$W_1(n)$ is isomorphic to a certain deformation of $W_1(1)\otimes \scr
O_1(n-1)$. Namely, $W_1(n)$ is isomorphic to the Lie algebra  $\germ L(\scr
O_1(n-1),-\partial )$. This isomorphism is given by $e_i\otimes x^k\mapsto
e_{pk+i}$ (the $e_i$'s on the left-hand side belong to $W_1(1)$, while the ones
on the right-hand side belong to $W_1(n))$.
\proclaim{Theorem 2.1} $H^2(\germ L(U,D))\cong \{G\in U^*|D\circ G=0\}$.
\par The basic cocycles may be chosen as follows{\rm:}
$$e_i\otimes u\wedge e_j\otimes v\mapsto
\delta_{i+j,0}(-1)^iG(vD(u))+\delta_{i+j,p}(-1)^iG(uv).\tag11$$
\endproclaim
\demo{Proof} It is a routine task to verify that cochains determined by
formula (11) are cocycles, and we leave it again to the reader. 
\par Let $\Phi\in Z^2(\germ L(U,D))$. Since $e_0\otimes \Bbb F$ is a torus in $
\germ L(U,D)$ and the corresponding root spaces are $\germ L_i=e_i\otimes U$, by
Lemma 1.1 we may assume that
$$\phi(e_i\otimes U,e_j\otimes U)=0\quad\text{if }i+j\not=0\text{ or
}p.\tag12$$
The cocycle equation is equivalent to the two following equalities:
$$\multline
\qquad\quad N_{i,j}\Phi(e_{i+j}\otimes uv,e_k\otimes w)\\
\qquad\qquad+(\text{simultaneous cyclic permutations 
of }(i,j,k)\\
\text{ and }(u,v,w))=0,\quad\{i,j,k\}\not=\{-1,-1,2\},\endmultline\tag13
$$
and
$$\aligned  
\Phi(e_{-1}\otimes& v,e_1\otimes uw)-\Phi(e_{-1}\otimes u,e_1\otimes
vw)\\
&-\Phi(e_2\otimes w,e_{p-2}\otimes vD(u)-uD(v))=0.\endaligned\tag14
$$
In (13), $k=-1$, $i+j=p+1$:
$$\Phi(e_i\otimes u, e_{j-1}\otimes vw)+\Phi(e_{i-1}\otimes uw, e_j\otimes
v)=0.$$
Substitute here $w=1$:
$$\Phi(e_i\otimes u,e_{j-1}\otimes v)=-\Phi(e_{i-1}\otimes u, e_j\otimes
v),\qquad i+j=p+1.\tag15$$
Substituting the last equality in the previous one, we obtain
$$\Phi(e_{i-1}\otimes u, e_j\otimes vw)=\Phi(e_{i-1}\otimes uw,e_j\otimes v).$$
Substitute here $v=1$:
$$\Phi(e_{i-1}\otimes u, e_j\otimes w)=\Phi(e_{i-1}\otimes uw, e_j\otimes
1),\qquad i+j=p+1.\tag16$$
(15) and (16) imply
$$\Phi(e_i\otimes u, e_j\otimes v)=(-1)^iG(uv),\qquad i+j=p,\tag17$$
for certain $G\in U^*$.
\par In (14), substitute $v=1$:
$$\Phi(e_{-1}\otimes u, e_1\otimes w)=\Phi(e_{-1}\otimes 1, e_1\otimes
uw)-G(wD(u)).\tag18$$
In (13), set $i=-1$, $j=0$, $k=1$, $w=1$:
$$\Phi(e_0\otimes u, e_0\otimes v)=\Phi(e_{-1}uv,e_1\otimes
1)-\Phi(e_{-1}\otimes u,e_1\otimes v).\tag19$$
Substituting the previous equality in the last one, we obtain
$\Phi(e_0\otimes u,e_0\otimes v)=-G(uD(v))$. Then the skew-symmetry of $\Phi$
implies
$-G(uD(v))-G(vD(u))=0$ which is equivalent to $D\circ G=0$.
\par Now put $\omega \in C^1(\germ L(U,D))$ as follows: $\omega (e_0\otimes
u)=\Phi(e_{-1}\otimes 1, e_1\otimes u)$; $\omega (e_i\otimes U)=0$, $i\not=0$.
Then
$d\omega (e_{-1}\otimes 1, e_1\otimes u)=\Phi(e_{-1}\otimes 1,e_1\otimes u)$.
Hence if we take $\Phi-d\omega $ instead of $\Phi$, we obtain a cocycle
satisfying all previous equalities and, moreover,
$$\Phi(e_{-1}\otimes 1, e_1\otimes u)=0.\tag20$$
\par Collecting (12) and (17)--(20), we obtain formula (11).
\par It remains to prove that if $\Phi=d\omega $, where $\Phi$ is determined by
formula (11) for an appropriate $G$, and $\omega \in C^1(\germ
 L(U,D))$, then
$G=0$. But this is obvious:
$$G(u)=\Phi(e_2\otimes u,e_{p-2}
\otimes 1)=d\omega (e_2\otimes u,e_{p-2}\otimes
1)=0.\qed$$
\enddemo
\proclaim{Corollary 2.2} $(p>3)\ H^2(W_1(\infty)\otimes U)=0$.
\endproclaim
\demo{Joint proof of Corollaries
 {\rm 1.5(i)} and \rm2.2} According to the
above-mentioned observation by Kuznetsov, $\germ L(U\otimes \scr
O_1(n-1),-1\otimes\partial )\cong W_1(n)\otimes U$, and the isomorphism
being given by $e_i\otimes u\otimes x^k\mapsto e_{pk+i}\otimes u$. 
\par Let $n<\infty$. By Theorem 2.1
$$\align
H^2(W_1(n)\otimes U)\cong&\{G\in(U\otimes \scr O_1(n))^*|G(U\otimes\langle
1,x,x^2,\dotsc, x^{p^{n-1}-2}\rangle )=0\}\\
\cong&U^*.\endalign
$$
Rewriting (11) in our case, we have
$$\align
&e_{pk+1}\otimes v\wedge e_{p1+j}\otimes u\mapsto-\pmatrix
k+l-1\\k-1\endpmatrix\delta_{i+j,0}(-1)^iG(uv\otimes x^{k+l-1})\\
&\qquad+\pmatrix k+l\\
k\endpmatrix \delta_{i+j,p}(-1)^iG(uv\otimes x^{k+l}),\qquad
{i,j=-1,0,1,\dotsc, p-2.}\endalign
$$
Direct easy calculations show that it is exactly the formula
$$e_i\otimes u\wedge e_j\otimes v\mapsto\delta_{i+j,p^n}(-1)^iG(uv).$$
\par Let $n=\infty$. Since $\partial (\scr O_1(\infty))=\scr O_1(\infty)$, by
Theorem 2.1 $H^2(W_1(\infty)\otimes U)=0$.\qed
\pagebreak
\enddemo
\remark{Remark} The expression of basic cocycles of $H^2(W_1(n)\otimes U)$  in
terms of the basis $\{e_i\}$ corresponds to the one for $H^2(W_1(n))$ found by
Dzhumadil'daev \cite6.
\endremark
\heading3. Central extensions of modular semisimple lie algebras\endheading
\par According to \cite3, Lie algebras of type $\germ g\otimes \scr
O_n+1\otimes \germ D$ and $W_1(n)\otimes
 \scr O_n+1\otimes\germ D$, where $\germ D$ is a subalgebra of
$\operatorname{Der}(\scr O_n)$ such that $\scr O_n$ does not contain $\scr
D$-invariant ideals, are one of the simplest but quite typical examples of 
modular semisimple Lie algebras. Preceding results allow us to determine their
central extensions.
\proclaim{Lemma 3.1} Let $L+\germ D$ be a semidirect product of Lie algebras,
where $L$ is an ideal with $[L,L]=L$. Then $H^2(L+\germ D)\cong H^2(L)^{\germ D}
+H^2(\germ D)$.\endproclaim
\demo{Proof} We may identify in an obvious manner $C^2(L+\germ D)$ with
$C^2(L)+\mathbreak
(L\otimes \germ D )^*+C^2(\germ D)$. Let us see how this identification
affects cocycles and coboundaries. Let $\Phi\in C^2(L+\germ D)$ and $\Phi=\Psi+
\Upsilon +\Xi$ be the corresponding decomposition. The cocycle equation $d\Phi=0$
is equivalent to the following relations:
$$\gather
d\Psi(x, y,z)=0,\\
\Upsilon([x,y],D)=\Psi([x,D],y)+\Psi(x,[y,D]),\\
\Upsilon([x,D_1],D_2)-\Upsilon([x,D_2],D_1)-\Upsilon(x,[D_1,D_2])=0,\\
d\Xi(D_1,D_2,D_3)=0,\endgather
$$
for all $x,y,z\in L$ and $D,D_1,D_2,D_3\in \germ D$. The second relation shows
that $\Psi\in Z^2(L)^{\germ D} := \{\Psi\in
Z^2(L)|\Psi([x,D],y)+\Psi(x,[y,D])\in B^2(L)$ $\forall D\in\germ D\}$ and, since
$[L,L]=L$, $\Upsilon$ is fully determined by $\Psi$. Direct verification shows
that the third relation follows from the two previous ones. So we have $Z^2(L+
\germ D)\cong Z^2(L)^{\germ D}+Z^2(\germ D)$. Obviously $B^2(L+\germ D)\cong
B^2(L)+B^2(\germ D)$. The last two isomorphisms imply the assertion of the
lemma.\qed
\enddemo
\par Corollaries 1.4, 1.5, and the preceding lemma immediately imply the
following proposition.
\proclaim{Proposition 3.2} $(\roman i)$ $H^2(\germ g\otimes U+1\otimes \germ
D)\cong HC^1(U)^{\germ D}+H^2(\germ D)$.
\par$(\roman{ii})$ $(p>3,\ n<\infty)$, $H^2(W_1(n)\otimes U+1\otimes \germ
D)\cong \{G\in U^*|D\circ G=0\ \forall D\in\germ D \} +H^2(\germ
D)$.\qed
\endproclaim
\par We conclude with the computation of the first-order cyclic cohomology of
the truncated polynomial algebra $\scr O_n$. We use multi-index notation: $
\Gamma_n := \{\alpha =(\alpha _1,\dotsc, \alpha _n)\in\Bbb Z^n|0\le\alpha
_i<p\}$, $x^\alpha =x^{\alpha _1}_1\dotsb x^{\alpha _n}_n$, $\varepsilon _i$
denotes 
$$(0,\dotsc, \underset{(i)}\to1,\dotsc,0).$$
\proclaim{Proposition 3.3} $\dim HC^1(\scr O_n)=(n-1)p^n+1$. Each cocycle is of
the form
$$x^\alpha \wedge x^\beta \mapsto \sum^n_{i=1}\alpha _i\lambda _i(\alpha +\beta
-\varepsilon _i)\tag21$$
where $\lambda _i(\alpha )\in\Bbb F$ are parameters satisfying the relations 
$$\sum^n_{i=1}\alpha _i\lambda _i(\alpha -\varepsilon _i)=0\tag22$$
for any $\alpha \in\Gamma_n$.%\longpage{13pc}
\pagebreak\endproclaim
\demo{Proof} From the computational point of view it is more convenient to
consider cyclic homology. Since $\scr O_n$ is commutative, $HC_1(\scr O
_n)\cong \Omega^1_{\scr O_n}/d\scr O_n$, where $\Omega^1_{\scr O_n}$ is the
module of K\"ahler differentials (cf.\ \cite{5, 9}). It is easy to see that
each element $udv$, $u, v\in\scr O_n$, may be represented as a sum of elements
of kind $x^\alpha \,dx^{\varepsilon _i}$. Moreover, we have
$$0=dx^\alpha =\sum^n_{i=1}\alpha _ix^{\alpha -\varepsilon _i}\,dx_i\tag23$$
for each $\alpha \in\Gamma_n$. So there are $np^n$ generators of the vector
space $HC_1(\scr O_n)$, and $p^n-1$ linear relations between them (if $\alpha
=0$, (23) gives a trivial equality). Hence $\dim HC_1(\scr O_n)=(n-1)p^n+1$.
\par Because of the duality isomorphism $HC^1(\scr O_n)\cong HC_1(\scr O_n)^*$,
the first assertion of the lemma holds.
\par Further, it is easy to verify that the cochains, defined by (21) ad (22),
are really cocycles. Since there are $np^n$ parameters and $p^n-1$ relations
between them, the dimension equality implies that all cocycles are taken into
account.\qed
\enddemo
\par Given a concrete subalgebra $\scr D$, it is easy to compute $HC^1(U)^{
\scr D}$. In the general case we are able to give a ``lower bound'' only:
\proclaim{Proposition 3.4} $\dim HC^1(\scr O_n)^{\operatorname{Der}(\scr
O_n)}=np^{n-1}$.
\par Each $\operatorname{Der}(\scr O_n)$-invariant cocycle is determined by
formula {\rm (21)}, where
$$\lambda _i(\alpha )=0\quad\text{if }\alpha _i\not=p-1.\tag24$$
\endproclaim
\demo{Proof} $\operatorname{Der}(\scr O_n)$ has as a basis $\{x^\gamma\partial
/\partial x_j|\ \gamma\in\Gamma_n,\ 1\le j\le n\}$. Writing the $x^\gamma 
\partial /\partial x_j$-invariance condition for the cocycle determined by
formula (21), we get
$$\sum^n_{i=1}(\alpha _i\beta _j-\alpha _j\beta _i)\lambda _i(\alpha +\beta +
\gamma-\varepsilon _i-\varepsilon _j)=0.\tag25$$
Setting here $\beta =\varepsilon _j$ and $\gamma=0$, we obtain $\alpha
_j\lambda _j(\alpha -\varepsilon _j)=0$, which is equivalent to (24). Obviously,
the last condition implies (22) and (25).\qed
\enddemo
\par In particular, $HC^1(\scr O_1)^{\scr D}=HC^1(\scr O_1)$ for any $\scr
D\subseteq \operatorname{Der}(\scr O _1)$.
\heading Acknowledgment\endheading

I am deeply grateful to A. S. Dzhumadil'daev for
friendly support and stimulating discussions.
\remark{Note added in proof} Recently A. Haddi (Comm.\ Algebra {\bf 20} (1992),
1145--1166) has obtained the general formula expressing
the homology group\linebreak$H_2(L\otimes A)$ in terms of $H_2(L)$
and $HC_1(A)$ for arbitrary Lie algebra $L$ in the case
of characteristic zero.
\endremark

\enddocument

\title{Erratum and addendum to ``Central extensions of current algebras''}
\author{Pasha Zusmanovich}
\address{Reykjav\'ik Academy, Iceland}
\email{pasha@akademia.is}

%    The 2010 edition of the Mathematics Subject Classification is
%    now available.  If you are citing a classification from the
%    new scheme, use the following input coding instead.

% pasha: if defining 2010 here, I get warning that the 2010 subject class is undefined;
% I guess the amsart.cls in not uptodate
%\subjclass[2010]{Primary 17B56; Secondary 13D03, 17B50}
\date{January 2, 2010; last revised April 25, 2010}

\maketitle

Regretfully, the author has spotted, with almost two decades delay, errors in \cite{tams}.

In \cite[Proposition 3.3]{tams}, the dimension of the first order cyclic cohomology 
$HC^1(O_n)$ of the reduced polynomial algebra 
$O_n=K[x_1, \dots, x_n]/(x_1^p, \dots, x_n^p)$ over a field $K$ of characteristic $p>0$, 
is specified incorrectly.

The correct value can be easily obtained as follows. 
For $n=1$, we have $\dim HC^1(O_1) = 1$ (see, for example,
\cite[Corollary 5.4.17]{loday}), the basic cocycle being
\begin{equation}\label{alpha}
\alpha(x^i, x^j) = \begin{cases}
i, &\text{if } i+j=p \\
0, &\text{otherwise} .
\end{cases}
\end{equation}
Now, from the fact that 
\begin{equation}\label{tensor}
O_n \simeq O_{n-1} \otimes O_1 ,
\end{equation}
and the K\"unneth formula for the cyclic cohomology (see formula (\ref{kunneth}) below), 
we get by induction $\dim HC^1(O_n) = np^{n-1}$.

\cite[Proposition 3.4]{tams} is wrong. The correct statement is given below. 
Recall that for an associative commutative algebra $A$, its derivation algebra $Der(A)$ 
acts on $HC^1(A)$  by the formula
\begin{equation*}
(D \varphi)(a,b) = \varphi(D(a),b) + \varphi(a,D(b))
\end{equation*}
for any $D \in Der(A)$, $\varphi\in HC^1(A)$, $a,b\in A$.

\begin{proposition}\label{1}
$$
HC^1(O_n)^{Der(O_n)} = 
\begin{cases}
HC^1(O_1), &n=1    \\
0,         &n>1 .
\end{cases}
$$
\end{proposition}

\begin{proof}
In the case $n=1$, it is easy to check that the cocycle (\ref{alpha}) is invariant under
the action of $x^k \frac{d}{dx}$ for any $0\le k < p$.

\begin{lemma}
Let $A,B$ be two associative commutative algebras with unit. Then
\begin{multline*}
HC^1(A\otimes B)^{Der(A\otimes B)} 
\\ \simeq 
\Big(HC^1(A)^{Der(A)} \otimes (B^*)^{Der(B)}\Big) 
\oplus
\Big((A^*)^{Der(A)} \otimes HC^1(B)^{Der(B)}\Big) .
\end{multline*}
\end{lemma}

\begin{proof}
By the K\"unneth exact sequence for the cyclic cohomology (see \cite[\S 4.4.10]{loday}), 
\begin{equation}\label{kunneth}
HC^1(A \otimes B) \simeq 
\Big(HC^1(A) \otimes B^*\Big) \oplus \Big(A^* \otimes HC^1(B)\Big) ,
\end{equation}
each cyclic $1$-cocycle on $A\otimes B$ being represented as the linear span of 
cocycles of the form
\begin{equation}\label{cocycle}
(a \otimes b) \wedge (a^\prime \otimes b^\prime) \mapsto
\varphi(a,a^\prime) \beta(bb^\prime) + \alpha(aa^\prime) \psi(b,b^\prime) ,
\end{equation}
where $\varphi\in HC^1(A)$, $\psi\in HC^1(B)$, $\alpha\in A^*$, $\beta\in B^*$, and
$a,a^\prime \in A$, $b,b^\prime \in B$.

Similarly, by the K\"unneth formula for the Hochschild 
cohomology (derivations are just Hochschild cocycles of order $1$),
$$
Der(A \otimes B) \simeq \Big(Der(A) \otimes B\Big) \oplus \Big(A \otimes Der(B)\Big) ,
$$
each derivation of $A\otimes B$ being represented as the linear span of derivations of 
the form $D \otimes R_b + R_a \otimes F$, where $D\in Der(A)$, $F\in Der(B)$, 
$a\in A$, $b\in B$, and $R_u$ denotes the multiplication by an element $u$ in the 
respective algebra.

Using these explicit isomorphisms, we get that the result of the action of the derivation 
$D \otimes 1_B$, where $D\in Der(A)$ and $1_B$ is the unit of $B$, on a cocycle
of the form (\ref{cocycle}), is equal to
$$
(D\varphi)(a,a^\prime) \otimes \beta(bb^\prime) + 
\alpha(D(aa^\prime)) \otimes \psi(b,b^\prime) ,
$$
and hence
$$
HC^1(A \otimes B)^{Der(A) \otimes 1_B} \simeq 
\Big(HC^1(A)^{Der(A)} \otimes B^*\Big) \oplus \Big((A^*)^{Der(A)} \otimes HC^1(B)\Big) .
$$
Similarly,
$$
HC^1(A \otimes B)^{1_A \otimes Der(B)} \simeq 
\Big(HC^1(A) \otimes (B^*)^{Der(B)}\Big) \oplus \Big(A^* \otimes HC^1(B)^{Der(B)}\Big) .
$$
Taking intersection of the right-hand sides in the last two formulas, we arrive at
the right-hand side of the isomorphism claimed. It is straightforward to check
that each cocycle of the form (\ref{cocycle}), where $\varphi\in HC^1(A)^{Der(A)}$,
$\psi\in HC^1(B)^{Der(B)}$, $\alpha\in (A^*)^{Der(A)}$, $\beta\in (B^*)^{Der(B)}$,
is invariant under the action of $Der(A\otimes B)$.
\end{proof}

\begin{remark}
Similar, but much more involved arguments could be used to prove the analogous
statements for not necessarily commutative algebras, for algebras without unit, and for 
the cyclic cohomology of higher degrees. As we are interested here solely
in application to the central extensions of modular semisimple Lie algebras, 
as specified below, we will not go into detail.
\end{remark}

\noindent\textit{Continuation of the proof of Proposition \ref{1}}. 
The case $n>1$ follows from (\ref{tensor}), the lemma just proved, and the obvious
fact that $(Der(O_n))(O_n) = O_n$ for any $n$.
\end{proof}

These errors do not affect other statements in \cite{tams}.

A more general, then those presented in \cite{tams}, formula for the second cohomology 
with trivial coefficients of the current Lie algebra, 
allows us to obtain results more general than \cite[Proposition 3.2]{tams}.
Namely, for modular semisimple Lie algebras which are represented as the semidirect
sum $(S \otimes O_n) \inplus \mathscr D$, where $S$ is a simple Lie algebra, and
$\mathscr D$ is a Lie subalgebra of $Der(O_n)$, we have:

\begin{proposition}
$$
H^2\Big((S \otimes O_n) \inplus \mathscr D\Big) \simeq
\Big(H^2(S) \otimes (O_n^*)^{\mathscr D}\Big) \oplus 
\Big(\mathcal B(S) \otimes HC^1(O_n)^{\mathscr D}\Big) \oplus H^2(\mathscr D) .
$$
\end{proposition}

Here $\mathcal B(S)$ denotes the space of symmetric invariant bilinear forms on $S$.

\begin{proof}
This is an immediate application of \cite[Lemma 3.1]{tams} (which, in its turn, is a simple consequence 
of the Hochschild--Serre spectral sequence) and a cohomological version of 
\cite[Theorem 0.1]{asterisque}.
\end{proof}

It is interesting to compare this situation with a somewhat dual one in \cite[\S 6]{dz}.
There, we have a similar formula for the space of so-called commutative $2$-cocycles 
on semisimple modular Lie algebras of the form $(S \otimes O_n) \inplus \mathscr D$, 
but we do not bother with the cyclic cohomology, as the corresponding direct summand
vanishes due to vanishing of its first tensor factor -- a certain skew-symmetric analog
of $\mathcal B$.